\documentclass{siamltex1213}

\usepackage{amsmath}
\usepackage{amssymb}

\newcommand{\bi}[1]{\textbf{\textit{#1}}}

\title{A 4th-Order Particle-in-Cell Method with Phase-Space Remapping for the Vlasov-Poisson Equation \thanks{This material is based upon work supported by the U.S. Department of Energy, Office of Science, Advanced Scientific Computing Research Program and performed under the auspices of the U.S. Department of Energy by Lawrence Berkeley National Laboratory under Contract DE-AC02-05CH11231. This work relied heavily on yt \cite{yt_paper} and the core scientific Python packages, including SciPy \cite{scipy} IPython \cite{Ipython}, NumPy \cite{NumPy}, and Matplotlib \cite{Matplotlib}, for data analysis and plotting.}}

\author{A. Myers \footnotemark[3] \footnotemark[4] \thanks{AM wishes to thank Anshu Dubey, Daniel Graves, and Daniel Martin for sharing their time and expertise with Chombo.} \and P. Colella  \footnotemark[3] \and B. Van Straalen  \footnotemark[3]}

\begin{document}
       
\maketitle
\slugger{sicomp}{xxxx}{xx}{x}{x--x}%slugger should be set to mms, siap, sicomp, sicon, sidma, sima, simax, sinum, siopt, sisc, or sirev
   
\renewcommand{\thefootnote}{\fnsymbol{footnote}}   
   
\footnotetext[3]{Applied Numerical Algorithms Group, Lawrence Berkeley National Laboratory, MS 50A-1148, Berkeley, CA, 94720}
\footnotetext[4]{atmyers@lbl.gov}

\renewcommand{\thefootnote}{\arabic{footnote}}
    
\begin{abstract}
Numerical solutions to the Vlasov-Poisson system of equations have important applications to both plasma physics and cosmology. In this paper, we present a new Particle-in-Cell (PIC) method for solving this system that is 4th-order accurate in both space and time. Our method is a high-order extension of one presented previously [B. Wang, G. Miller, and P. Colella, SIAM J. Sci. Comput., 33 (2011), pp. 3509--3537]. It treats all of the stages of the standard PIC update - charge deposition, force interpolation, the field solve, and the particle push - with 4th-order accuracy, and includes a 6th-order accurate phase-space remapping step for controlling particle noise. We demonstrate the convergence of our method on a series of one- and two- dimensional electrostatic plasma test problems, comparing its accuracy to that of a 2nd-order method. As expected, the 4th-order method can achieve comparable accuracy to the 2nd-order method with many fewer resolution elements.
\end{abstract}


    
    \section{Introduction}
        
     In this paper, we present a method for solving the Vlasov-Poisson system of equations, which in non-dimensional form is given by:
\begin{equation}
\label{eq:vlasov}
\frac{\partial f}{\partial t} + \bi{v}\cdot  \frac{\partial f}{\partial \bi{x}} - \bi{E} \cdot \frac{\partial f}{\partial \bi{v}} = 0
\end{equation}
and
\begin{equation}
\label{eq:poisson}
\nabla^2 \phi = -\rho.
\end{equation} 

Here, $f(\bi{x}, \bi{v}, t)$ is the phase space distribution function defined on $(\bi{x}, \bi{v}) \in \mathbb{R}^D \times \mathbb{R}^D$, where $D = 1, 2,$ or $3$ is the number of spatial dimensions under consideration, $\bi{E}(\bi{x}, t) = - \nabla  \phi$ is the electric field, $\phi(\bi{x}, t)$ is the potential, and $\rho(\bi{x}, t)$ is the charge density. For simplicity, we assume that the our system contains a single, negatively charged species, and that the ions form a fixed, neutralizing background, so that the charge density can be defined as: 
\begin{equation}
\label{eq:density}
\rho(\bi{x}, t) = 1 - \int_{\mathbb{R}^D} f(\bi{x}, \bi{v}, t) d \bi{v}.
\end{equation}
     
    Equations (\ref{eq:vlasov}) and (\ref{eq:poisson}) describe the phase space evolution of a collisionless fluid under the influence of electrostatic forces. They are important for plasma physics, where they are used, for example, to model space plasmas, particle accelerators, and for controlled thermonuclear fusion. In a slightly different form, Equations (\ref{eq:vlasov}) and (\ref{eq:poisson}) can also be applied to cosmology, where they are used to model the gravitational evolution of dark matter in an expanding universe. For simplicity, we have specialized to the plasma version of the Vlasov-Poisson system in this paper, but the methods presented here can be easily applied to the self-gravitating case as well.
    
 The Vlasov Equation (\ref{eq:vlasov}) is a nonlinear advection equation in phase space and can in principle be solved with a variety of grid-based methods \cite{filbet_comparison_2003, banks_new_2010}, including high-order methods \cite[e.g.]{vogman_dory_2014}. However, particle discretizations, which reduce the Vlasov-Poisson system to a set of coupled ordinary differential equations, have been more common in practice. The Particle-in-Cell (PIC) method \cite{hockney_computer_1981}, in which the forces are computed on an intermediate grid and then interpolated back to the particle positions, is a particularly simple approach that has been widely used in both cosmology \cite[e.g.]{heitmann_robustness_2005} and plasma physics \cite[e.g.]{warp}. 
 
Traditional PIC methods, however, suffer from a few downsides. The first is that they are usually limited to 2nd-order accuracy. While 4th-order PIC methods have been developed in the context of fluid simulation \cite{edwards_2012}, and high-order field solves and time integrators have been used in the plasma context \cite[e.g.]{jacobs_2006}, to our knowledge, ours is the first PIC method for Vlasov-Poisson that treats \emph{every} stage of the PIC update with fourth-order accuracy. 

The second downside concerns the stability of PIC methods over long time evolutions. The convergence theory for electrostatic PIC \cite{wang_particle--cell_2011} shows that the stability error for the electric field contains a term that grows exponentially with time. While the growth rate of this term is problem-dependent, given enough time it can eventually degrade the accuracy of the solution, a problem often described as ``particle noise." The PIC method presented in \cite{wang_particle--cell_2011}, however, circumvents this problem by periodically restarting the calculation with a new set of particles that represent the same underlying distribution function. Such \emph{remapping} or \emph{remeshing} techniques have also been applied successfully in the context of fluid dynamics to vortex methods \cite{cottet_vortex_2000} and to smoothed particle hydrodynamics \cite{chaniotis_remeshed_2002}. With remapping, PIC methods can obtain accurate numerical solutions to the Vlasov-Poisson problem for long time integrations in both the plasma \cite{wang_particle--cell_2011, wang_adaptive_2012} and the cosmological \cite{myers2015} context.

However, the PIC method in \cite{wang_particle--cell_2011} and \cite{wang_adaptive_2012} was only 2nd-order accurate. For obtaining accurate numerical solutions with a feasible number of resolution elements, higher-order methods are greatly desirable. This is particularly true given current trends in high-performance computing. As computer architectures evolve, the limiting factor affecting application performance is increasingly not the rate at which the processor can perform computations, but rather than rate at which data can be streamed to the processor from DRAM. In light of this, the Arithmetic Intensity (AI) - the number of arithmetic operations per byte of DRAM accessed - of a given numerical algorithm is considered to be a critical factor in achieving the theoretical maximum performance on current and next-generation supercomputing platforms \cite{williams_adaptive_2009}. 2nd-order PIC algorithms have peak arithmetic intensities of around 1 Flop/byte, putting them in the streaming-limited regime of current and planned supercomputing architectures and making it impossible to achieve maximum performance. High-order methods, however, which allow for more operations per byte, can potentially achieve peak performance on these machines. 

In this paper, we extend the PIC algorithm of \cite{wang_particle--cell_2011} and \cite{wang_adaptive_2012} to be 4th-order accurate in both space and time. The heart of our method is the use of the high-order interpolating functions derived in \cite{lo_real_2015} for charge deposition and force interpolation. We have implemented this algorithm in Chombo \cite{chombo_design} and used it to solve a suite of test problems, demonstrating its fourth-order convergence rates. The outline of the paper is as follows. We begin by giving a review of PIC methods in Section \ref{sec:overview}. We then describe our 4th-order PIC method, as well as a 2nd-order PIC method that we use for comparison, in Sections \ref{sec:2nd_order_method} and \ref{sec:4th_order_method}. We present our numerical results on set of one- and two-dimensional test problems in Section \ref{sec:results}. Finally, in Section \ref{sec:conclusions}, we present our conclusions and discuss future research directions.
    
\section{Numerical Methods}
\label{sec:methods}
    
\subsection{An overview of PIC methods}
\label{sec:overview}

To solve equations (\ref{eq:vlasov}) and (\ref{eq:poisson}) with PIC, the initial distribution function must be sampled by a set of Lagrangian particles, $\mathbb{P}$. In this paper, we generate this initial particle distribution using the approach outlined in \cite{wang_particle--cell_2011}. The particles are initially laid out at the cell centers of a Cartesian grid in phase space with mesh spacings $h_x$ and $h_v$ in position and velocity space, respectively. For a given initial distribution function, $f(\bi{x}, \bi{v}, t = 0)$, the initial particle representation is computed by assigning each particle $p \in \mathbb{P}$ a charge $q_p$:
\begin{equation}
q_p = f(\bi{x}_p^i, \bi{v}_p^i, t=0) h_x^D h_v^D.
\end{equation} where $\bi{x}_p^i$ and $\bi{v}_p^i$ are the initial position and velocity of particle $p$. For computational efficiency, we discard particles with charges less than some problem-dependent threshold value. The initial distribution function can then be approximated as:
\begin{equation}
    \label{discrete_f}
f(\bi{x}, \bi{v}, t = 0) \approx \sum_{p \in \mathbb{P}} q_p \delta \left(\bi{x} - \bi{x}_p^i \right) \delta \left(\bi{v} - \bi{v}_p^i \right).
\end{equation}

Once the particles have been generated, the equations of motion for their trajectories \(\left( \bi{x}_p(t), \bi{v}_p(t) \right)\) can be obtained by substituting Equation (\ref{discrete_f}) into Equation (\ref{eq:vlasov}). The result is the following system of ODEs:
\begin{align}
\label{eq:motion}
\frac{d q_p}{dt} &= 0 \nonumber \\
\frac{d \bi{x}_p}{dt} &= \bi{v}_p \nonumber \\
\frac{d \bi{v}_p}{dt} &=  -\bi{E}_p,
\end{align}
where $\bi{E}_p$ is the acceleration on particle $p$ induced by the surrounding charge distribution. 

Equation $\ref{eq:motion}$, when coupled with a procedure for computing the electric field at the particle positions, can then be used to numerically advance Equations (\ref{eq:vlasov}) and (\ref{eq:poisson}) in time. In PIC methods, this time advance is computed in a number of stages, as follows. First, the charge density is computed on a grid via a deposition step, in which each particle distributes its charge into a number of cells using some interpolating function. These functions are commonly taken to be one of the $B$-splines, e.g. the ``Cloud-In-Cell" (CIC) and ``Triangle-Shaped Cloud" (TSC) functions. However, while higher-order $B$-splines have increasing degrees of smoothness, they are limited to interpolating with at most 2nd-order accuracy \cite{monaghan_particle_1985}. Next, Poisson's equation for the potential is solved on the grid using a finite difference method along with some kind of fast Poisson solver, such as fast Fourier Transforms or multigrid methods. Next, the electric field is computed on the grid using a finite difference approximation to the gradient. Once the field is known at the grid points, it can be interpolated back to the particle positions. It is important that this be done using the same interpolating function as in the deposition step, in order to avoid spurious self-forces. Finally, once the electric field at the particle positions is known, the particles can be advanced in time using some numerical ODE solver. If this solver contains multiple stages, such as with Runge-Kutta methods, this entire procedure will need to be completed several times per time step.

The final ingredient needed for accurate PIC calculations, as discussed in \cite{wang_particle--cell_2011} and \cite{wang_adaptive_2012}, is a particle remapping procedure. During the remap, the current set of particles is replaced by a new set that encodes the same underlying distribution function. This process is similar to the particle initialization procedure described above, except that instead of generating the particles from a given initial distribution function, we generate them by depositing the known particle distribution onto a high-dimensional, Cartesian grid in phase space and sampling the resulting distribution at each cell center. It is important to note that, unlike with pure grid methods, this process does not require storing the entire high-dimensional grid in phase space at once; rather, because each particle only contributes charge to nearby positions, it can be done on purely local chunks in phase space. It is also important to note that this process does not need to be done every time step. Once the new particles have been generated, the calculation can continue. This procedure prevents errors in the particle trajectories from compounding to the point that they significantly degrade the solution.

They key factors affecting the accuracy of PIC methods are thus 1) the accuracy of the interpolating function used for charge deposition and force interpolation, 2) the accuracy of the finite difference stencils used for the field solve, 3) the accuracy of the time integration scheme, and 4) the accuracy of the interpolating function used for the remap. In what follows, we first review a basic PIC method that is 2nd-order accurate (very similar to the one from \cite{wang_particle--cell_2011}), and then we describe a new PIC method that extends the first scheme to 4th-order accuracy.

\subsection{A 2nd-order PIC method}
\label{sec:2nd_order_method}

The stages of the 2nd-order PIC algorithm proceed as follows.

\begin{itemize}

\item \textbf{Particle Deposition} The charge density is defined on a Cartesian mesh $ \bi{x}_{\bi{i}} = \left( \bi{i} + 1/2 \right) \Delta \bi{x}$, 
where $ \bi{i} \in \mathbb{Z}^D $ are the cell indices and $\Delta \bi{x} $ is the cell spacing. At second order, we can use the following deposition step:
\begin{equation}
\rho_{\bi{i}} = \sum_p \left( \frac{q_p}{V_i} \right) \bi{W}_{\bf{2}} \left( \frac{\bi{x}_{i} - \bi{x}_p}{\Delta x} \right).
\end{equation} Here, $V_i = {\Delta x}^D$ is the volume of cell $i$ and $\bi{W}_{\bf{2}}(\bi{x})$ is a $D$-dimensional, ``Cloud-in-Cell" interpolating function, given by: 

\begin{equation}
\bi{W}_{\bf{2}} \left(\bi{x} \right) = \prod_{d = 1}^D W_{2} \left( x_d \right),
\end{equation}

\begin{equation}
   W_{2}(x) = \left\{
     \begin{array}{lr}
       1 - \lvert x \rvert, & 0 \leq \lvert x \rvert \leq 1, \\
       0 & \text{otherwise.}
     \end{array}
   \right.
\end{equation}

Note that in general, we do not use the same mesh spacing for the particle discretization and the Poisson mesh, i.e. $\Delta x \ne h_x$.

\item \textbf{Field Solve} The next step is to solve the Poisson equation for the potential at the same grid points on which the density is defined. At 2nd order, we discretize the Laplacian operator using the standard \(2 D\)+\(1\) point centered difference approximation:
\begin{equation}
- \sum_{d=1}^D \frac{\phi_{\bi{i} + \bi{e}^d} - 2\phi_{\bi{i}} + \phi_{\bi{i} - \bi{e}^d}}{{\Delta  x}^2} = -\rho_{\bi{i}}.
\end{equation}

We solve the resulting linear system with geometric multigrid, using Gauss-Seidel with Red-Black ordering as the smoother. Through this paper, we use periodic boundary conditions and set the solver tolerance to $10^{-9}$. After iterating to convergence, the electric field components can be computed on the mesh as
\begin{equation}
\bi{E}_{\bi{i}}^{d} = - \frac{\phi_{\bi{i} + \bi{e}^d} - \phi_{\bi{i} - \bi{e}^d}}{2 \Delta x}.
\end{equation}

\item \textbf{Force Interpolation} Next, we interpolate the field back to the particle positions using the same interpolating function as in the deposition step:

\begin{equation}
\bi{E}_p = \sum_{\bi{i}} \bi{E}_{\bi{i}} V_i \bi{W}_{\bf{2}} \left( \frac{\bi{x}_i -\bi{x}_p}{\Delta x} \right).
\end{equation}

\item \textbf{Particle Push} The final stage is the update of the particle positions. To numerically integrate Equation (\ref{eq:motion}), we use the following 2nd-order accurate Runge-Kutta (RK2) method:
\begin{eqnarray}
\bi{x}^{n+1} &=& \bi{x}^n + \frac{1}{2} \bi{k}_1 \Delta t^2 \nonumber \\
\bi{v}^{n+1} &=& \bi{v}^n + \frac{1}{2} \left( \bi{k}_1 + \bi{k}_2 \right) \Delta t,
\end{eqnarray} where

\begin{eqnarray}
\bi{k}_1 &=& \bi{E}(\bi{x}^n) \nonumber \\
\bi{k}_2 &=& \bi{E}(\bi{x}^n + \bi{v}^n \Delta t).
\end{eqnarray} For each stage of the RK2 time step, we must generate a right-hand side for the Poisson equation corresponding to the appropriate particle positions by performing particle deposition, the field solve, and the force interpolation steps. 
    
 \item \textbf{Remapping} The remapping step involves replacing the current set of particles with a new set that is laid out on a Cartesian grid in phase space. As in particle initialization, the cell spacings of this mesh are $h_x$ and $h_v$ in position and velocity space, and we only keep particles with weights that exceed some problem-dependent threshold value. The charges of the new particles $q_p^{*}$ are then calculated as 
\begin{equation}
q_p^{*} = \sum_p q_p \bi{W}_3 \left( \frac{\bi{x}_\bi{i} - \bi{x}_p}{h_x} \right) \bi{W}_3 \left( \frac{\bi{v}_{\bi{j}} - \bi{v}_p}{h_v} \right),
\end{equation} where $\bi{W}_3(\bi{x})$ is the 3rd-order interpolating function from \cite{wang_particle--cell_2011} and \cite{wang_adaptive_2012}, given by
\begin{equation}
\bi{W}_3 \left(\bi{x} \right) = \prod_{d = 1}^{D} W_{3} \left( x_d \right),
\end{equation} and
\begin{equation}
   W_3(x) = \left\{
     \begin{array}{lr}
       1 - \frac{5}{2} \lvert x \rvert^2 + \frac{3}{2} \lvert x \rvert^3, & 0 \leq \lvert x \rvert \leq 1, \\
       \frac{1}{2}\left( 2 - \lvert x \rvert \right)^2 \left( 1 - \lvert x \rvert \right), & 1 \leq \lvert x \rvert \leq 2, \\ 
       0 & \text{otherwise.}
     \end{array}
   \right.
\end{equation} A 3rd-order interpolant is necessary here in order for the overall scheme to be 2nd-order accurate, since one order of accuracy is lost in the remap step \cite{wang_particle--cell_2011}.
    
\end{itemize}
    
    \subsection{A 4th-order PIC method}
    \label{sec:4th_order_method}
    
    We now describe our new, 4th-order method. This time, the PIC update proceeds as follows:

\begin{itemize}

\item \textbf{Particle Deposition} At fourth order, we can no longer use a $B$-spline to interpolate from the particle positions to the mesh cells. Instead, we use the following deposition step:

\begin{equation}
\rho_{\bi{i}} = \sum_p \left( \frac{q_p}{V_i} \right) \bi{W}_{\bf{4}} \left( \frac{\bi{x}_{i} - \bi{x}_p}{\Delta x} \right),
\end{equation} where $\bi{W}_{\bf{4}}(\bi{x})$ is the $D$-dimensional, 4th-order accurate interpolating function from \cite{lo_real_2015}, given by: 

\begin{equation}
\bi{W}_{\bf{4}} \left(\bi{x} \right) = \prod_{d = 1}^D W_{4} \left( x_d \right),
\end{equation}

\begin{equation}
   W_{4}(x) = \left\{
     \begin{array}{lr}
       1 - \frac{\lvert x \rvert}{2} - {\lvert x \rvert}^2 + \frac{{\lvert x \rvert}^3}{2}, & 0 \leq \lvert x \rvert \leq 1,\nonumber \\
       1 - \frac{11 {\lvert x \rvert}}{6} + {\lvert x \rvert}^2 - \frac{{\lvert x \rvert}^3}{6}, & 1 \leq \lvert x \rvert \leq 2, \nonumber \\
       0 & \text{otherwise.}
     \end{array}
   \right.
\end{equation}

\item \textbf{Field Solve} As before, we solve the Poisson equation for the electrostatic potential at the same grid points on which the density is defined. This time, we discretize the Laplacian operator using a 4th-order centered-difference approximation:

\begin{equation}
-\sum_{d=1}^D \frac{- \phi_{\bi{i} + 2 \bi{e}^d} + 16 \phi_{\bi{i} + \bi{e}^d} - 30\phi_{\bi{i}} + 16 \phi_{\bi{i} - \bi{e}^d} - \phi_{\bi{i} + 2 \bi{e}^d}}{{12 \Delta  x}^2} = \rho_{\bi{i}},
\end{equation} which we solve using geometric multigrid, as before. The electric field can also be computed using 4th-order centered differences as:
\begin{equation}
\bi{E}_{\bi{i}}^d = - \frac{-\phi_{\bi{i} + 2 \bi{e}^d} + 8 \phi_{\bi{i} + \bi{e}^d} - 8 \phi_{\bi{i} - \bi{e}^d} + \phi_{\bi{i} + 2 \bi{e}^d}}{12 \Delta x}.
\end{equation}
\item \textbf{Force Interpolation} The force is computed at the particle positions as:

\begin{equation}
\bi{E}_p = \sum_{\bi{i}} \bi{E}_{\bi{i}} V_i \bi{W}_{\bf{4}} \left( \frac{\bi{x}_i -\bi{x}_p}{\Delta x} \right).
\end{equation}

\item \textbf{Particle Push} Instead of RK2, we use a 4th-order accurate Runge-Kutta (RK4) method. 
This method assumes the special case of a velocity independent force, and thus has only 3 stages instead of the usual four \cite{abramowitz_handbook_1972}:

\begin{eqnarray}
\bi{x}^{n+1} &=& \bi{x}^n + \bi{v}^n \Delta t + \frac{1}{6} \left(\bi{k}_1 + 2 \bi{k}_2 \right) \Delta t^2 \nonumber \\
\bi{v}^{n+1} &=& \bi{v}^n + \frac{1}{6} \left( \bi{k}_1 + 4 \bi{k}_2 + \bi{k}_3 \right) \Delta t,
\end{eqnarray} where

\begin{eqnarray}
\bi{k}_1 &=& \bi{E}(\bi{x}^n) \nonumber \\
\bi{k}_2 &=& \bi{E}(\bi{x}^n + \frac{1}{2} \bi{v}^n \Delta t + \frac{1}{8} \bi{k}_1 \Delta t^2) \nonumber \\
\bi{k}_3 &=& \bi{E}(\bi{x}^n + \bi{v}^n \Delta t + \frac{1}{2} \bi{k}_2 \Delta t^2).
\end{eqnarray} As before, for each of the three stages of the RK4 time step, we must generate a right-hand side for the Poisson equation corresponding to the appropriate particle positions by performing the particle deposition, field solve, and force interpolation steps. 

\item \textbf{Remapping}
    
As before, in order to preserve the overall accuracy of the method, the remapping step must use a spatial interpolation method that is at least one order higher than the desired order of the method as a whole. Thus, to apply remapping in our 4th-order PIC method, we use the following 6th-order interpolation function from \cite{lo_real_2015} to perform the remap:

\begin{equation}
\bi{W}_{\bf{6}} \left(\bi{x} \right) = \prod_{d = 1}^D W_{6} \left( x_d \right),
\end{equation}

\begin{equation}
   \label{eq:W6}
   W_{6}(x) = \left\{
     \begin{array}{lr}
       1 - \frac{{\lvert x \rvert}}{3} - \frac{5 {\lvert x \rvert}^2}{4} + \frac{5 {\lvert x \rvert}^3}{12} + \frac{{\lvert x \rvert}^4}{4} - \frac{{\lvert x \rvert}^5}{12}, & 0 \leq \lvert x \rvert \leq 1,\nonumber \\
       1 - \frac{13 {\lvert x \rvert}}{12} - \frac{5 {\lvert x \rvert}^2}{8} + \frac{25 {\lvert x \rvert}^3}{24} - \frac{ 3 {\lvert x \rvert}^4}{8} + \frac{{\lvert x \rvert}^5}{24}, & 1 \leq \lvert x \rvert \leq 2, \nonumber \\
       1 - \frac{137 {\lvert x \rvert}}{60} + \frac{15 {\lvert x \rvert}^2}{8} - \frac{17 {\lvert x \rvert}^3}{24} + \frac{{\lvert x \rvert}^4}{8} - \frac{{\lvert x \rvert}^5}{120}, & 2 \leq \lvert x \rvert \leq 3, \nonumber \\
       0 & \text{otherwise.}
     \end{array}
   \right.
\end{equation}

The new particle charges are then
\begin{equation}
q_p^{*} = \sum_p q_p \bi{W}_6 \left( \frac{\bi{x}_\bi{i} - \bi{x}_p}{h_x} \right) \bi{W}_6 \left( \frac{\bi{v}_{\bi{j}} - \bi{v}_p}{h_v} \right),
\end{equation}

\end{itemize} 

\subsection{Positivity Preservation}

\begin{figure}[h]
  \centering
    \includegraphics[width=\textwidth]{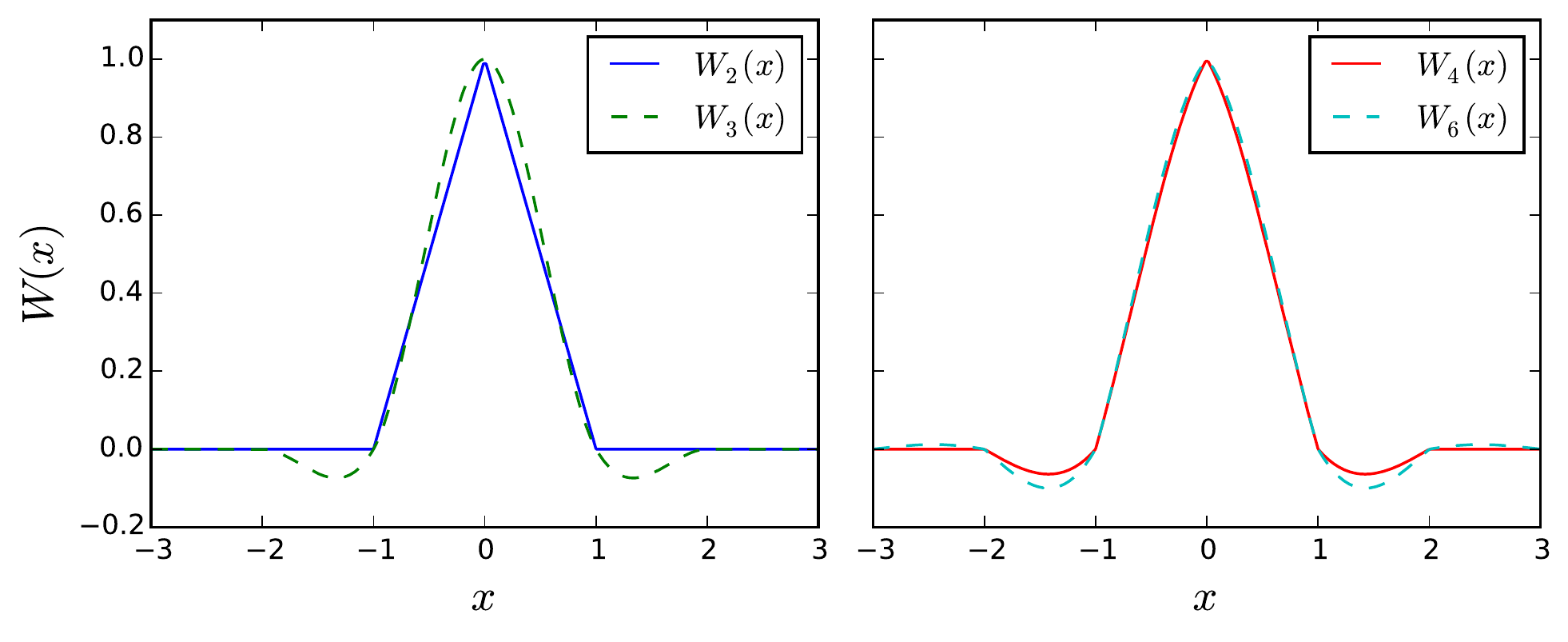}
    \caption{The interpolating kernels used in the various deposition / interpolation operations in this paper. Left panel - the kernels used by the 2nd-order method. Right panel - the kernels used by the 4th-order method. \label{fig:discrete_deltas}}
\end{figure}

We display all of the various interpolating functions we use in this work in Figure \ref{fig:discrete_deltas}. It is important to note that, unlike 2nd-order interpolants, higher-order interpolating functions are not strictly positivity preserving, in that a single particle with negative charge will not produce a uniformly negative charge density when deposited onto a grid. This can happen during the remap stage of the 2nd-order method, and during both the charge deposition and remap stages of the 4th-order method. One way to account for this is to apply a mass-redistribution algorithm to the distribution function during the remap step, following \cite{wang_particle--cell_2011}. That is, once we have generated the distribution function on a phase-space patch, we redistribute the undershoot in cell $\bi{i}$,
\begin{equation}
\delta f_{\bi{i}} = \min(0, f_{\bi{i}}^n)
\end{equation} to its neighbors $\bi{i} + \bi{l}$ in proportion to their capacity $\xi$, 
\begin{equation}
\xi_{\bi{i}+\bi{l}} = \max(0, f_{\bi{i}+\bi{l}}^n) 
\end{equation} so that
\begin{equation}
f_{\bi{i}+\bi{l}}^{n+1} = f_{\bi{i}+\bi{l}}^{n} + \frac{\xi_{\bi{i}+\bi{l}}}{\sum_{\bi{k} \ne 0}^{\text{neighbors}} \xi_{\bi{i}+\bi{k}}} \delta f_{\bi{i}}.
\end{equation} Here, $n$ and $n+1$ refer to the value before and after redistribution. In general, this procedure might need to be repeated for several iterations before the distribution function is strictly positive. In the problems we consider in this paper, we find that 2 or 3 iterations is sufficient.

We perform this redistribution procedure during the remapping phase on all the runs presented in Section \ref{sec:results}. However, we find that turning this off makes only a negligible difference on the test problems considered below. Because this may not be true for all problems, however, particularly those with large gradients in the distribution function, we include the positivity preservation algorithm here for completeness.
 
\subsection{Implementation}

We have implemented both of the above algorithms using the Chombo \cite{chombo_design} software framework for solving partial differential equations, which includes tools for doing parallel particle simulations that were originally developed in \cite{miniati_block_2007}. The parallelization uses standard message passing with block-structured domain decomposition in physical space. This code, as well as our data analysis scripts, is available online \footnote{https://bitbucket.org/atmyers/4thOrderPIC}. Our one-dimensional results presented in this paper were generated on a Macbook Pro, while the two-dimensional results were run on NERSC's Edison machine using up to 256 MPI processes.

\section{Numerical Results}
\label{sec:results}

In this section, we compare the performance of the 2nd- and 4th-order PIC methods on a set of standard test problems. We begin by considering two one-dimensional test problems, and then move on to perform similar tests in two spatial dimensions. In these numerical tests, we have fixed the ratio of $\Delta x / h_x$ to be 2. To compute convergence rates, we have used Richardson extrapolation to compute our error estimates. That is, if $\bi{E}^h$ is the electric field computed at a given spatial resolution and time step, and $\bi{E}^{2 h}$ is the electric field computed with \emph{all} the discrete elements of the problem ($\Delta x$, $h_x$, $h_v$, and $\Delta t$) coarsened by a factor of 2, then the solution error is defined as
\begin{equation}
e^h = \lvert \bi{E}^h - \bi{E}^{2 h} \rvert,
\end{equation} and the order of the method $q$ is computed as 
\begin{equation}
q = \log_2 \left( \frac{\lvert \lvert e^{2 h} \rvert \rvert}{\lvert \lvert e^{h} \rvert \rvert} \right). 
\end{equation} When spatial interpolation is required to compare the electric field between runs with different resolutions, we have used cubic spline interpolation (as implemented in SciPy \cite{scipy}), so as not to mask the 4th-order convergence results.

\subsection{1D Linear Landau Damping}
\label{sec:landau}

The first test problem we consider is Linear Landau Damping - the damped propagation of a small-amplitude plasma wave. To begin, we perform the calculation in one spatial dimension. We take the initial distribution function to be
\begin{align}
f(x, v, t = 0) &= \frac{1}{\sqrt{2 \pi}} \exp{\left(-v^2/2 \right)} \left( 1 + \alpha \cos \left( k x \right) \right) \nonumber \\
(x, v) &= \left[0, L = 2 \pi / k \right] \times \left[ -v_{\text{max}}, v_{\text{max}} \right],
\end{align}
where the amplitude of the perturbation $\alpha$ is $0.01$, its wavenumber $k$ is $0.5$, and $v_{\text{max}} = 10$. This problem uses periodic boundary conditions on the physical space domain $x \in (0, L)$. During particle initialization and during each remap, we have discarded particles with weights less than $10^{-16}$.

\begin{figure}[h]
  \centering
    \includegraphics[width=\textwidth]{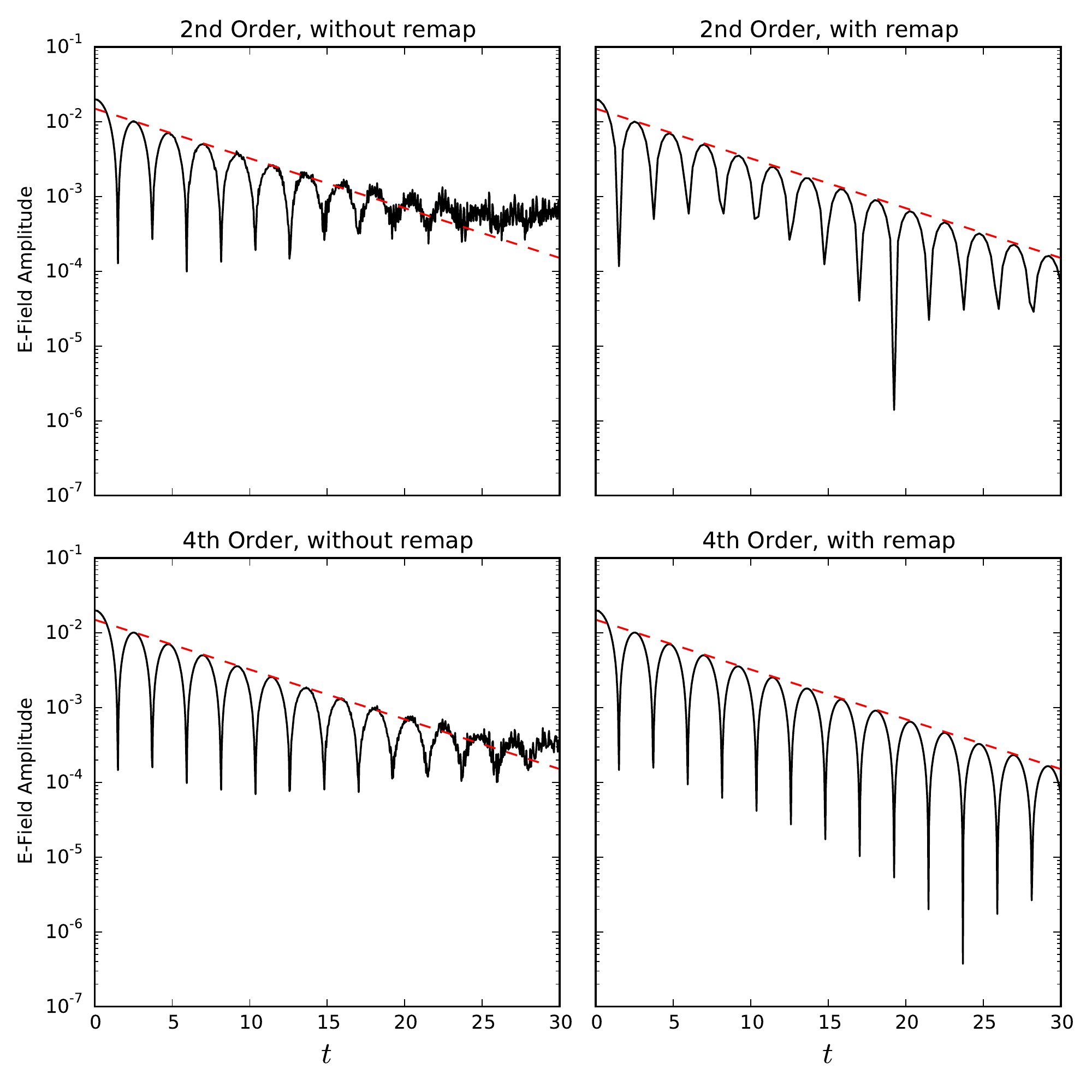}
    \caption{Amplitude of the electric field as a function of time for the one-dimensional linear Landau damping problem. The black solid line shows the numerical solution, while the red dotted line shows the theoretical damping rate. Top left - 2nd-order method, without remapping. Top Right - 2nd-order method, with remapping. Bottom Left - 4th-order method, without remapping. Bottom Right - 4th-order method, with remapping. All numerical calculations were performed with $N_{\text{cells}} = 64$, $N_{x} = 128$, and $N_{v} = 256$, and $dt = 1 / 32$. \label{fig:landau_solution}}
\end{figure}

This problem has been used extensively as a test  for plasma PIC codes. According to the analytic theory, the electric field is supposed to oscillate with a frequency $\omega = 1.416$, and the amplitude is supposed to be exponentially damped at a rate $\gamma = 0.1533$. We compare the results of our 2nd and 4th order PIC methods, both with and without remapping, to the analytic theory in Figure \ref{fig:landau_solution}. The numerical solutions used $N_{\text{cells}} = 64$ cells for the Poisson solve, $N_x = 128$ and $N_v = 256$ for the initial particle grid, and a PIC time step of $dt = 1 / 32$. In the runs that used remapping, we applied the remap every 5 time steps. All of the runs track the expected damping rate well at early times. As in  \cite{wang_particle--cell_2011}, however, the runs without remapping become noisy and fail to track the exact damping rate at late times. For both the 2nd-order and the 4th-order PIC methods, remapping appears to be necessary for long time evolutions.

\begin{figure}[h]
  \centering
    \includegraphics[width=\textwidth]{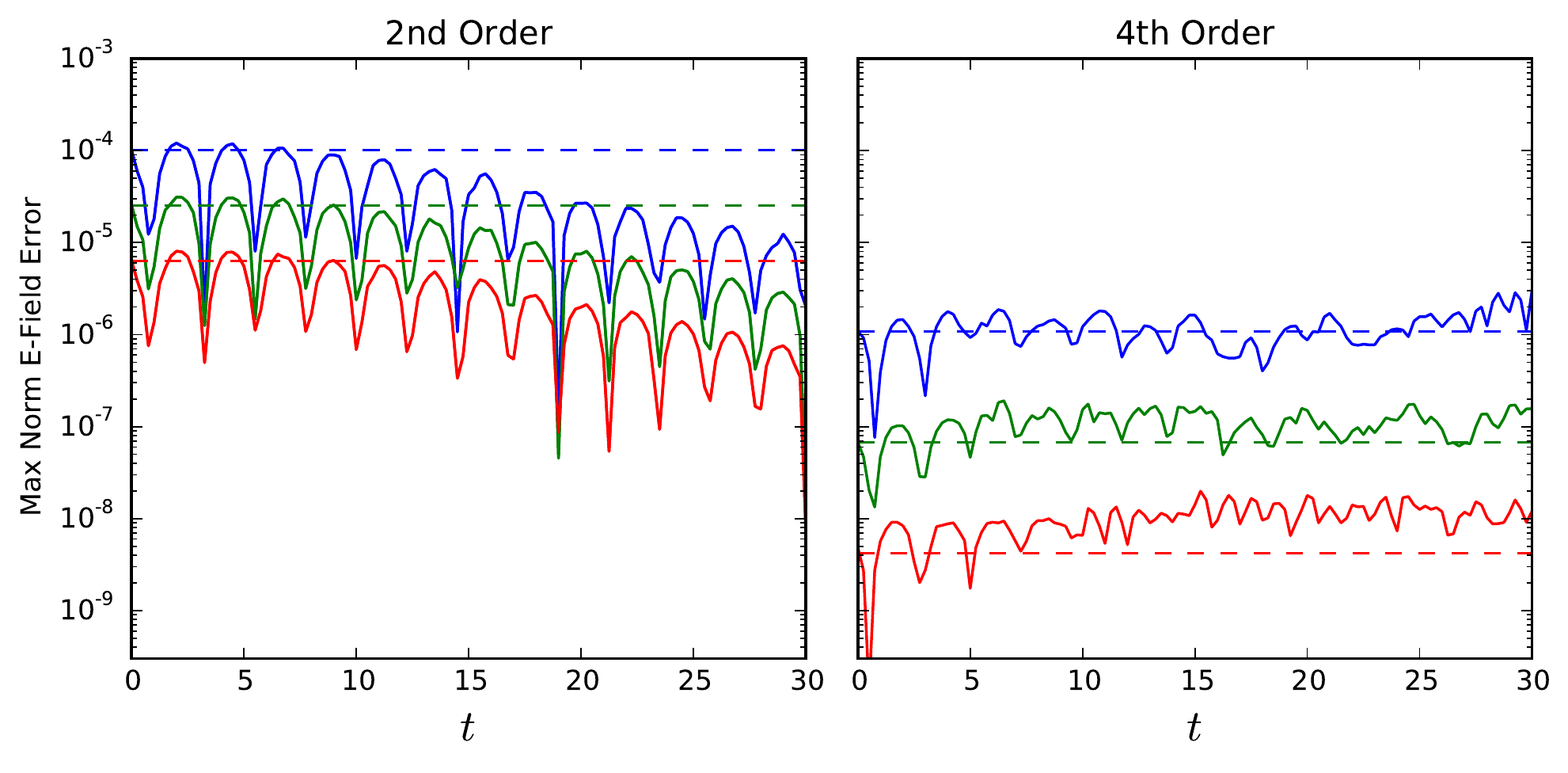}
    \caption{Max norm of the error in the electric field as a function of time for the one-dimensional linear Landau damping problem. Left - 2nd order method. Right - 4th order method. In both plots, the solid colored lines refer to the Richardson errors associated with consecutive pairs of resolutions, progressing from low (blue), to middle (green), to high (red). The dotted colored lines show how the initial errors should decrease if the methods performed at exactly 2nd (left) and 4th (right) orders. \label{fig:landau_error}}
\end{figure}

Next, we compare the accuracy of the 4th-order and 2nd-order PIC methods on this problem. We conduct a resolution study starting at  $N_{\text{cells}} = 32$, $N_x = 64$, $N_v = 128$ and $dt = 1 / 16$. We do four runs total, doubling the resolution and halving the time step with each successive calculation. We plot the max norm of the error in the electric field versus time in Figure \ref{fig:landau_error}, where we have used Richardson extrapolation to the estimate the error between each consecutive pair of resolutions. As expected, the 4th-order method is much more accurate than the 2nd-order method when run at the same resolution, by as much as two orders of magnitude for the highest resolution studied. Alternatively, the lowest resolution 4th-order runs are about as accurate as the highest resolution 2nd-order runs on this test problem.  

\begin{figure}[h]
  \centering
    \includegraphics[width=\textwidth]{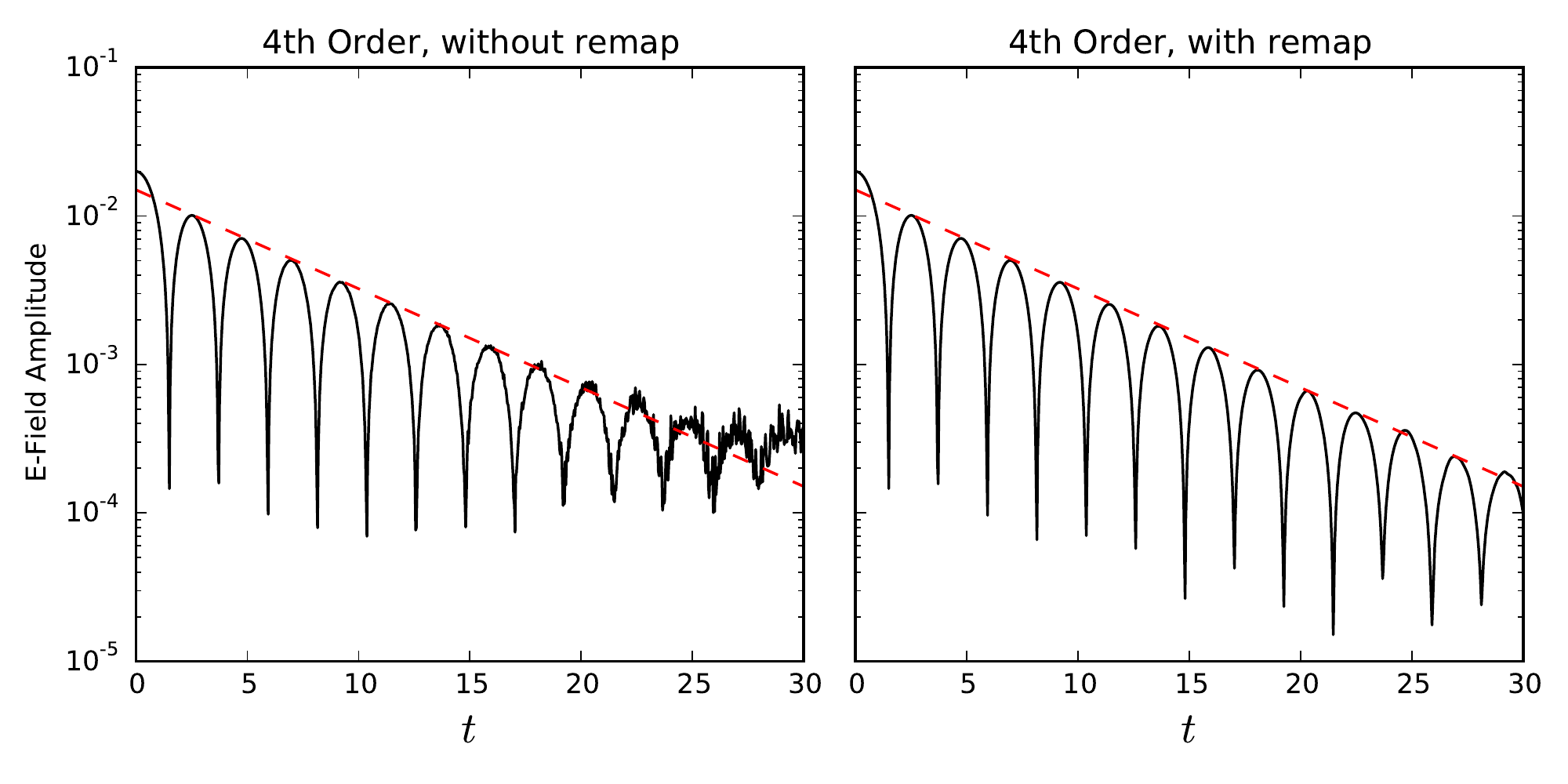}
    \caption{Amplitude of the electric field as a function of time for the one-dimensional linear Landau damping problem. The black solid line shows the numerical solution, while the red dotted line shows the theoretical damping rate. Left - 4th-order method, without remapping. Right - 4th-order method, remapping 9 times over the course of the simulation. All numerical calculations were performed with $N_{\text{cells}} = 64$, $N_{x} = 128$, and $N_{v} = 256$, and $dt = 1 / 32$. \label{fig:landau_solution_remap_compare}}
\end{figure}

The question of how often to apply the remapping procedure is an important consideration, and in general the answer will be problem dependent. In this test, we are concerned with demonstrating the 4th-order accuracy of our method. It is therefore crucial that the particle trajectory errors addressed by the remap be kept small compared to all the other sources of numerical error, so that 4th-order convergence can be observed. Experimentation suggests that, in order for this to happen, the remap needs to be applied every 5-10 time steps, which translates to remapping 100-200 times on the $N_{\text{cells}} = 256$ version of this test. However, for practical applications, it may not be necessary to keep the trajectory errors so small. Indeed, we find that remapping the particle distribution as few as 9 times over the course of the calculation greatly improves the situation, allowing us to track the expected damping rate all the way to $t = 30$  (see Figure \ref{fig:landau_solution_remap_compare}). 

We have enforced positivity preservation on all of our runs for this test problem; however, in this case, we find that turning it off makes no difference to the electric field to within machine precision.

\subsection{1D Two-Stream Instability}
\label{sec:twostream}
Next, we study the two-stream instability, again working in one spatial dimension. In this problem, there is a counter-streaming plasma flow in velocity space, along with a small initial density perturbation. The initial distribution function for this test is:

\begin{align}
f(x, v, t = 0) &= \frac{1}{\sqrt{2 \pi}} v^2 \exp{\left(-v^2/2 \right)} \left( 1 + \alpha \cos \left( k x \right) \right) \nonumber \\
(x, v) &= \left[0, L = 2 \pi / k \right] \times \left[ -v_{\text{max}}, v_{\text{max}} \right],
\end{align}

where we take $\alpha = 0.01$, $k = 0.5$. We have again adopted periodic boundary conditions in physical space and discarded particles with weights less than $10^{-16}$. 

\begin{figure}[h]
  \centering
    \includegraphics[width=\textwidth]{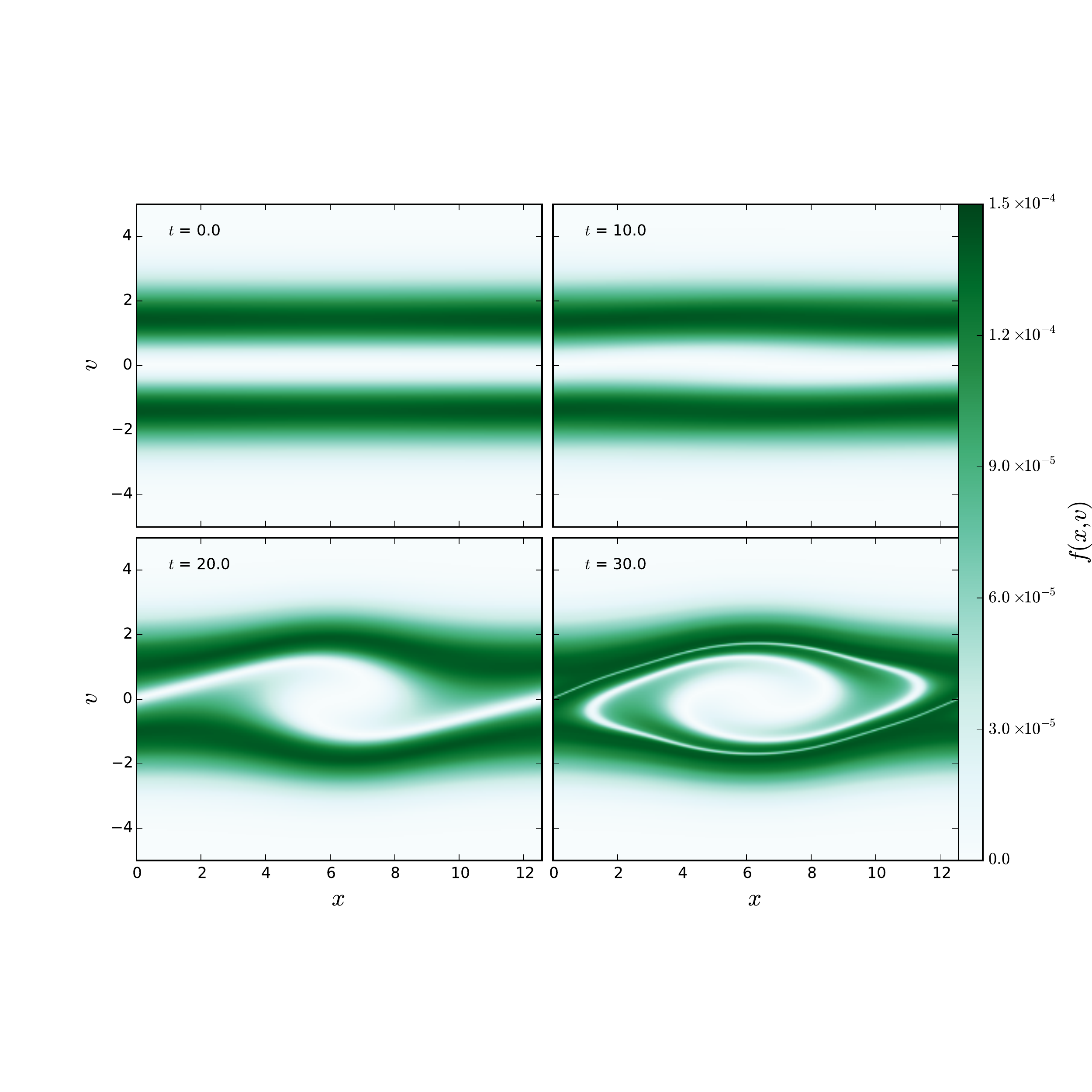}
    \caption{Phase-space distribution function at four selected times for the one-dimensional, two-stream instability problem. The data for this figure comes from the $N_{\text{cells}} = 256$, 4th-order run. For details of how this figure was generated, see Section \ref{sec:twostream} \label{fig:twostream_solution}}
\end{figure}

\begin{figure}[h]
  \centering
    \includegraphics[width=\textwidth]{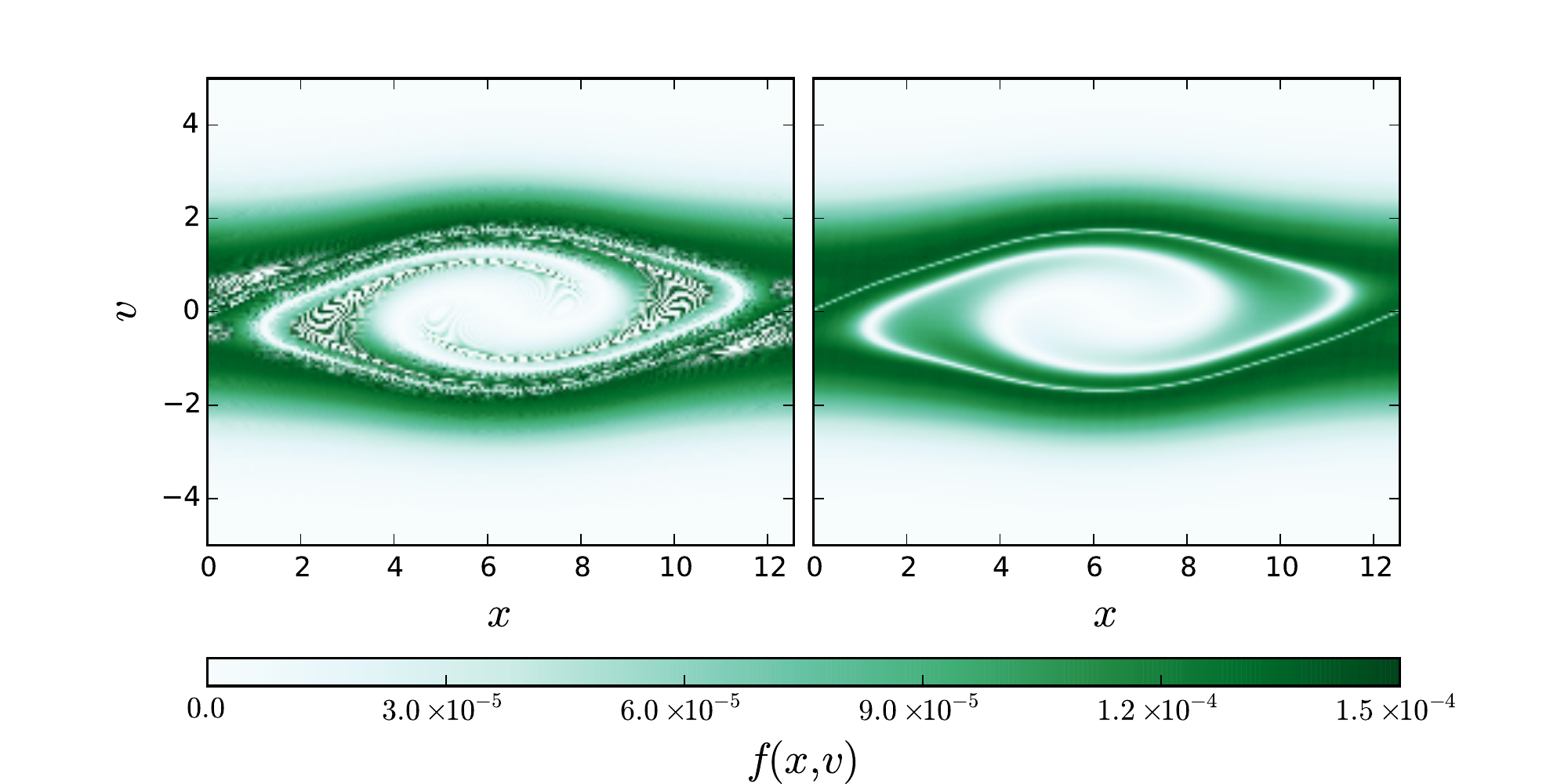}
    \caption{Phase-space distribution functions for the one-dimensional, two-stream instability problem. Left - 4th-order method, without remapping. Right - 4th-order method, remapping 9 times over the course of the simulation. \label{fig:twostream_solution_remap_compare}}
\end{figure}

To begin, we conduct a 4th-order run with $N_{\text{cells}} = 256$, $N_x = 512$, $N_v = 1024$, and $dt = 1 / 128$. The time evolution of the phase-space distribution function is shown in Figure \ref{fig:twostream_solution}. To construct the distribution function, we have used Equation (\ref{eq:W6}) to deposit the particles onto a $512$ by $1024$ mesh in phase space. We again apply the particle remap every 5 PIC time steps; however, as in Section \ref{sec:landau}, we find that applying the remap as few as 9 times over the simulation greatly reduces the degree of particle noise visible in the solution (Figure \ref{fig:twostream_solution_remap_compare}).

\begin{figure}[h]
  \centering
    \includegraphics[width=\textwidth]{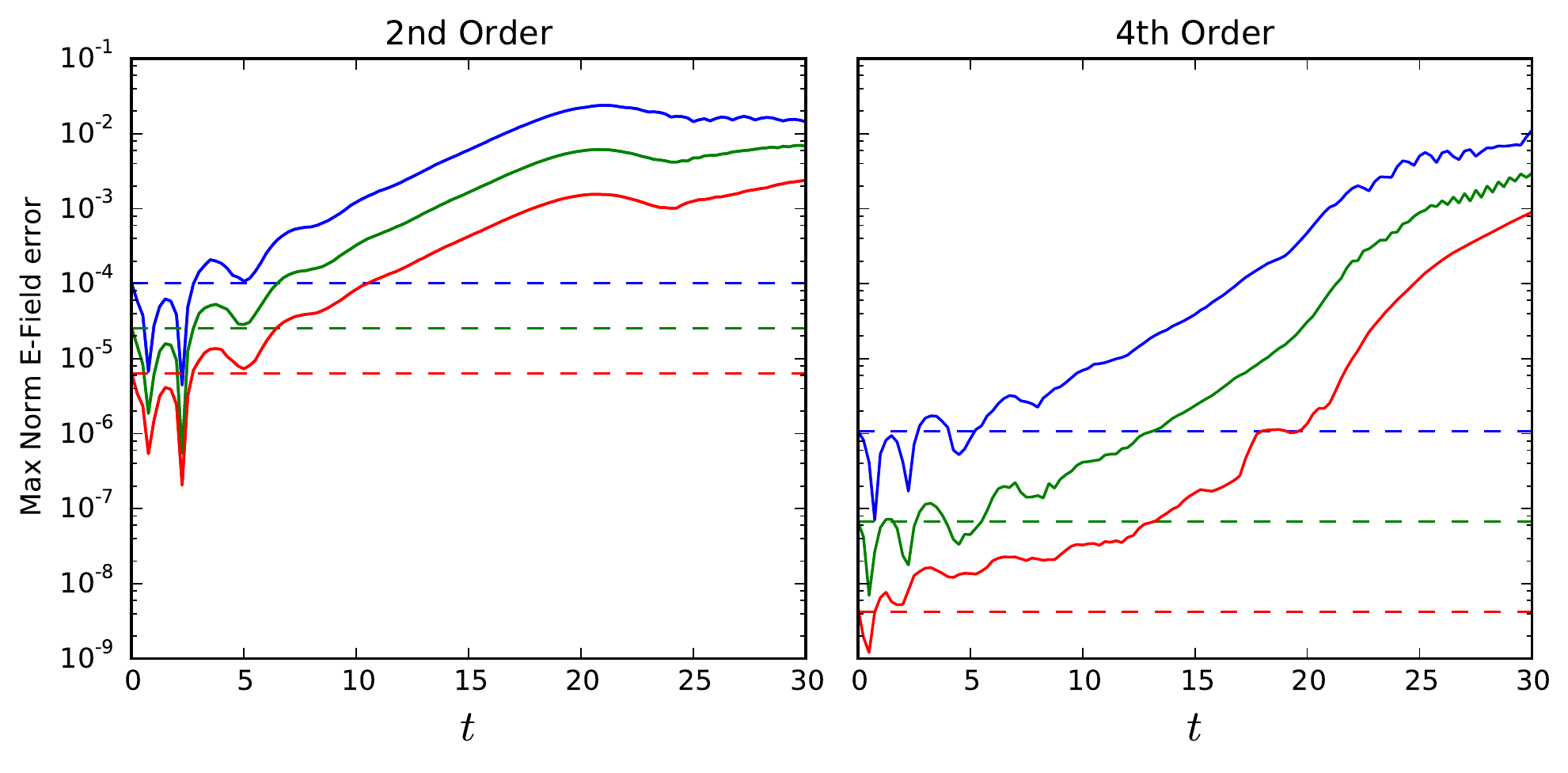}
     \caption{Max norm of the error in the electric field as a function of time for the one-dimensional, two-stream instability problem. Left - 2nd order method. Right - 4th order method. The line styles and colors have the same meaning as in Figure \ref{fig:landau_error}. \label{fig:twostream_error}}
\end{figure}

Next, we conduct a resolution study as in Section \ref{sec:landau}. We start at $N_{\text{cells}} = 32$, $N_x = 64$, $N_v = 128$, and $dt = 1 / 16$, and once again conduct four runs, doubling the resolution and decreasing the time step with each run. The resulting Richardson errors are shown in Figure \ref{fig:twostream_error}. As in the Landau damping example, at early times the 4th-order method is much more accurate than the 2nd-order method, by as much as two orders of magnitude at the resolutions we investigate. In the two-stream test, however, the errors grow significantly with time in both the 2nd-order and the 4th-order cases. This is due to the formation of thin, filamentary structures in phase space, as visible in the bottom right panel of Figure \ref{fig:twostream_solution}. These features are not well-resolved at the resolutions we use in our convergence study, and hence reduce the accuracy of the simulation from the expected order, such that, at very late times, the errors made by the two methods are comparable. However, the 4th-order method delays this loss of accuracy due to filamentation considerably. As late as time 20, for instance, the 4th order method is still around 2 orders more accurate than the 2nd-order method.

As before, we have included positivity preservation during the remap step on this test problem. Unlike with the linear Landau damping test, however, the change in the electric field is significantly higher than machine precision (as high as $ \sim 10^{-6}$ at late times) on this problem, particularly at late times when filamentation gives rise to large gradients in the distribution function. This difference is still small, however, compared to the other numerical errors (Figure \ref{fig:twostream_error}), and therefore running without positivity preservation does not affect the 4th-order convergence rates on this problem.

\subsection{2D Linear Landau Damping}
\label{sec:landau_2D}

The previous test problems both used one spatial dimension and one velocity dimension. Since certain numerical instabilities may only occur in higher-dimensional problems, it is also important to also check our algorithm on two-dimensional problems. To do so, we first repeat the linear Landau damping problem, this time performing the calculation in two-dimensional space. We take the initial distribution function to be
\begin{align}
f(x, y, v_x, v_y, t = 0) &= \frac{1}{2 \pi} \exp{\left(-(v_x^2 + v_y^2)/2 \right)} \left( 1 + \alpha \cos \left( k_x x \right) \cos \left( k_y y \right) \right) \nonumber \\
(x, y) &= \left[0, L = 2 \pi / k_x \right] \times  \left[0, L = 2 \pi / k_y \right] \nonumber \\
(v_x, v_y) &= \left[ -v_{\text{max}}, v_{\text{max}} \right] \times \left[ -v_{\text{max}}, v_{\text{max}} \right],
\end{align} where we set $\alpha = 0.05$, $k_x = k_y = 0.5$, $v_{\text{max}} = 6.0$, and use periodic boundary conditions in physical space.

\begin{figure}[h]
  \centering
    \includegraphics[width=\textwidth]{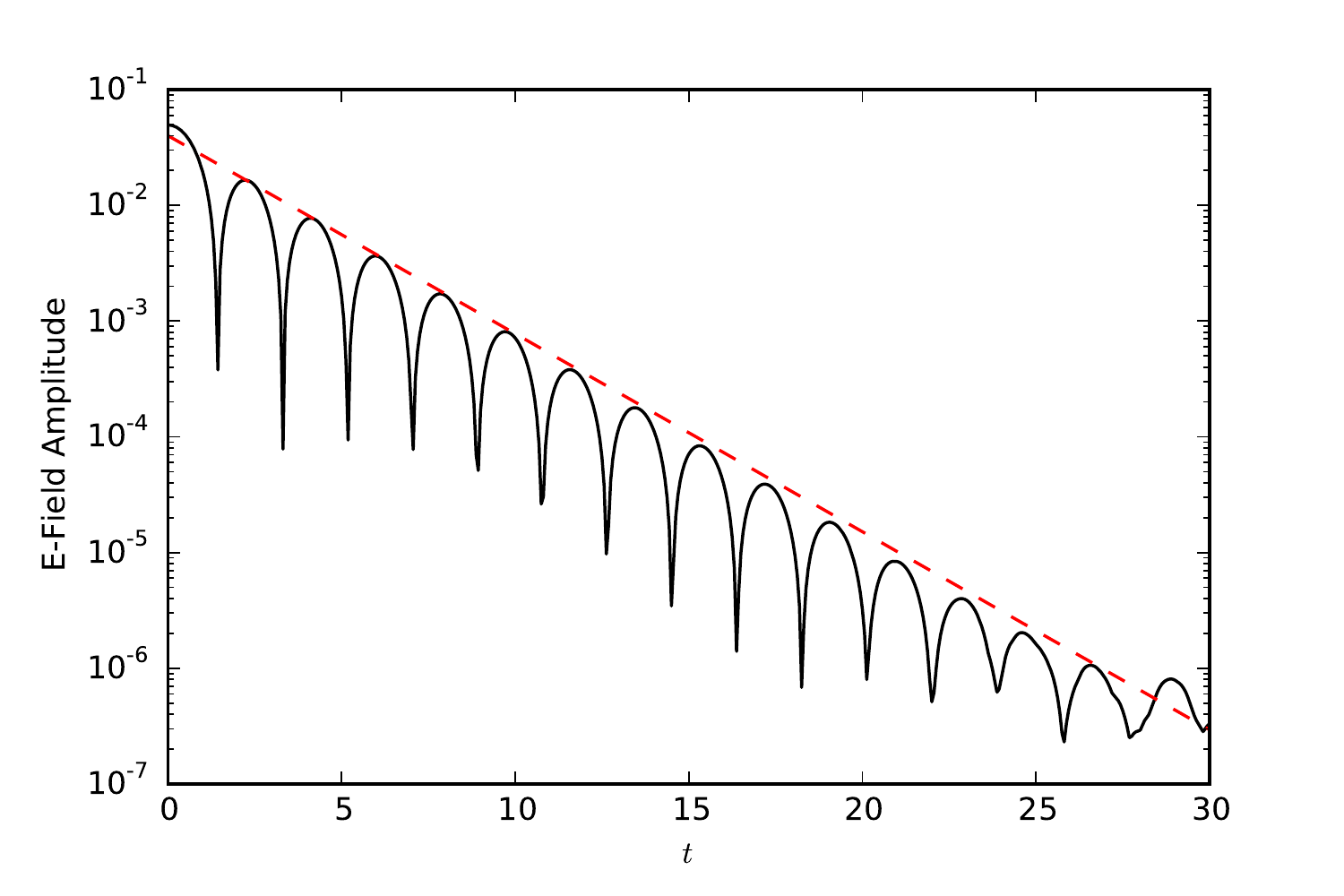}
    \caption{The amplitude of the electric field versus time for the 2D linear Landau damping problem. The black solid line shows the numerical result for the electric field, while the red-dotted line shows the expected damping rate of $\gamma = -0.394$. \label{fig:landau_solution_2D}}
\end{figure}

As before, this problem has an analytic solution for the damping rate of the electric field amplitude, which for these parameter choices should be $\gamma = -0.394$. To compare against this theory, we set $N_{\text{cells}} = 32$, $N_x = 64$, and $N_v = 128$, and $dt = 1/16$, and evolve the initial distribution out to time $t = 30.0$. We again apply remapping every 5 PIC time steps and discard particles with weights less than $10^{-12}$. The simulation electric field amplitude is compared against the theoretical expectation in Figure \ref{fig:landau_solution_2D}. Even at this relatively low resolution, we find that the numerical result track the expected damping rate well out to time $t \approx 28$. 

\subsection{2D Two-Stream Instability}
\label{sec:twostream_2D}

Finally, we perform the two-dimensional version of the two-stream instability problem.The initial distribution function is
\begin{align}
f(x, y, v_x, v_y, t = 0) &= \frac{1}{12 \pi} \exp{\left(-(v_x^2 + v_y^2)/2 \right)} \left( 1 + \alpha \cos \left( k_x x \right) \right) \left( 1 + 5 v_x^2\right) \nonumber \\
(x, y) &= \left[0, L = 2 \pi / k_x \right] \times  \left[0, L = 2 \pi / k_y \right] \nonumber \\
(v_x, v_y) &= \left[ -v_{\text{max}}, v_{\text{max}} \right] \times \left[ -v_{\text{max}}, v_{\text{max}} \right],
\end{align} where $\alpha = 0.05$, $k_x = k_y = 0.5$, $v_{\text{max}} = 9.0$, and we again use periodic boundary conditions in physical space.

\begin{figure}[h]
  \centering
    \includegraphics[width=\textwidth]{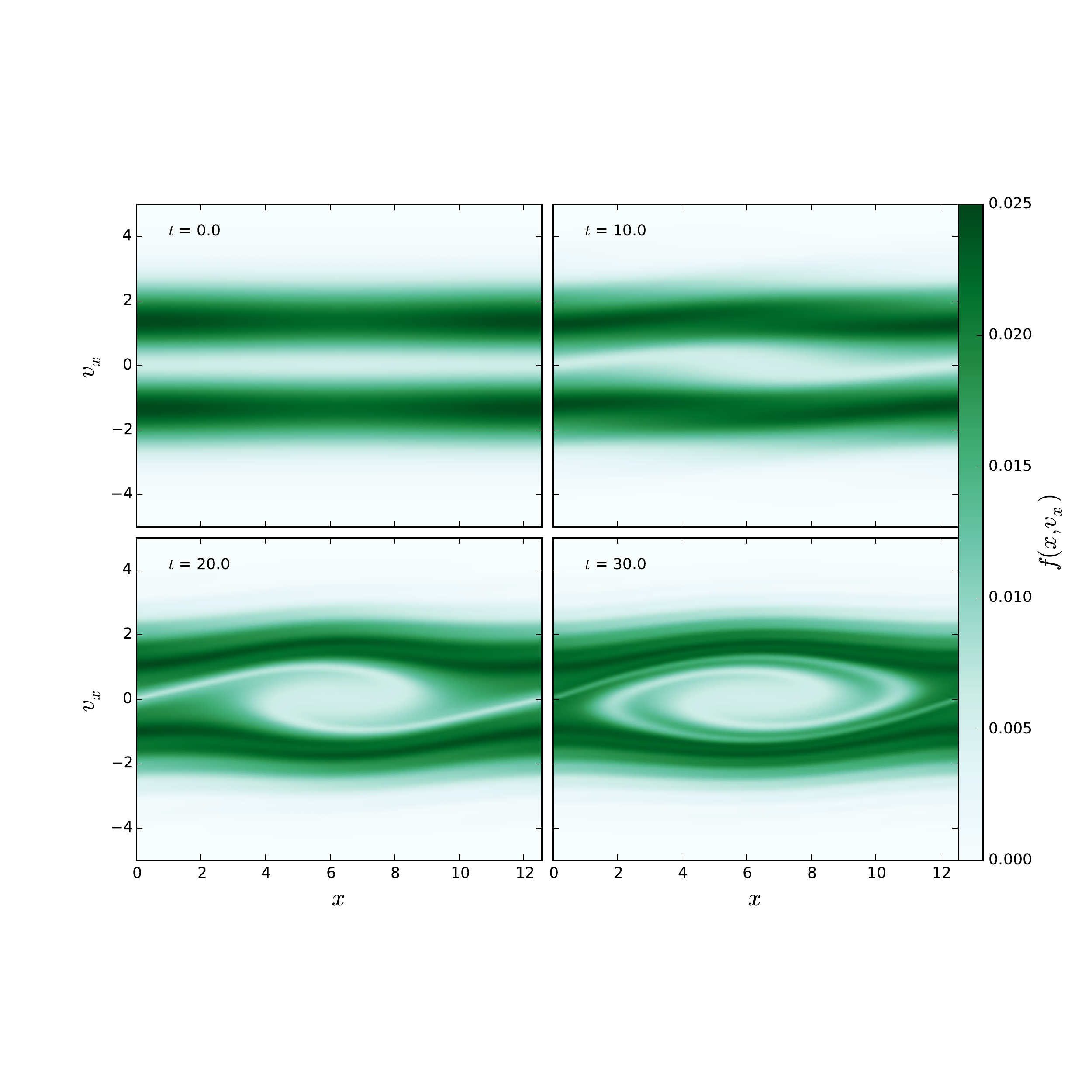}
    \caption{A 2D representation of the 4D phase space distribution function at four select times for the two-dimensional two-stream instability problem. \label{fig:twostream_solution_2D}} 
\end{figure}

We set $N_{\text{cells}} = 64$, $N_x = 128$, and $N_v = 256$, and $dt = 1/32$, and evolve the initial conditions out to time $t = 30.0$. As before, we apply remapping every 5 PIC time steps and discard particles with weights less than $10^{-12}$. The time evolution of a two-dimensional slice through the phase-space distribution function is shown in Figure \ref{fig:twostream_solution_2D}. The generate the 2$D$ version of the distribution in $(x, v_x)$ space, we deposit the particles using Equation (\ref{eq:W6}) onto a 2$D$ mesh with 128 by 256 grid points. That is, we compute
\begin{equation}
f_{i, j} = \sum_p \left( \frac{q_p}{\Delta x \Delta v} \right) W_6 \left( \frac{x_i - x_p}{\Delta x} \right) W_6 \left( \frac{v_{x,i} - v_{x,p}}{\Delta v} \right).
\end{equation} This is equivalent to depositing the particles onto a 4$D$ mesh and integrating over $y$ and $v_y$. As before, when remapping is employed, our 4th-order method tracks the evolution of the distribution function without visible particle noise. 

\begin{figure}
  \centering
    \includegraphics[width=\textwidth]{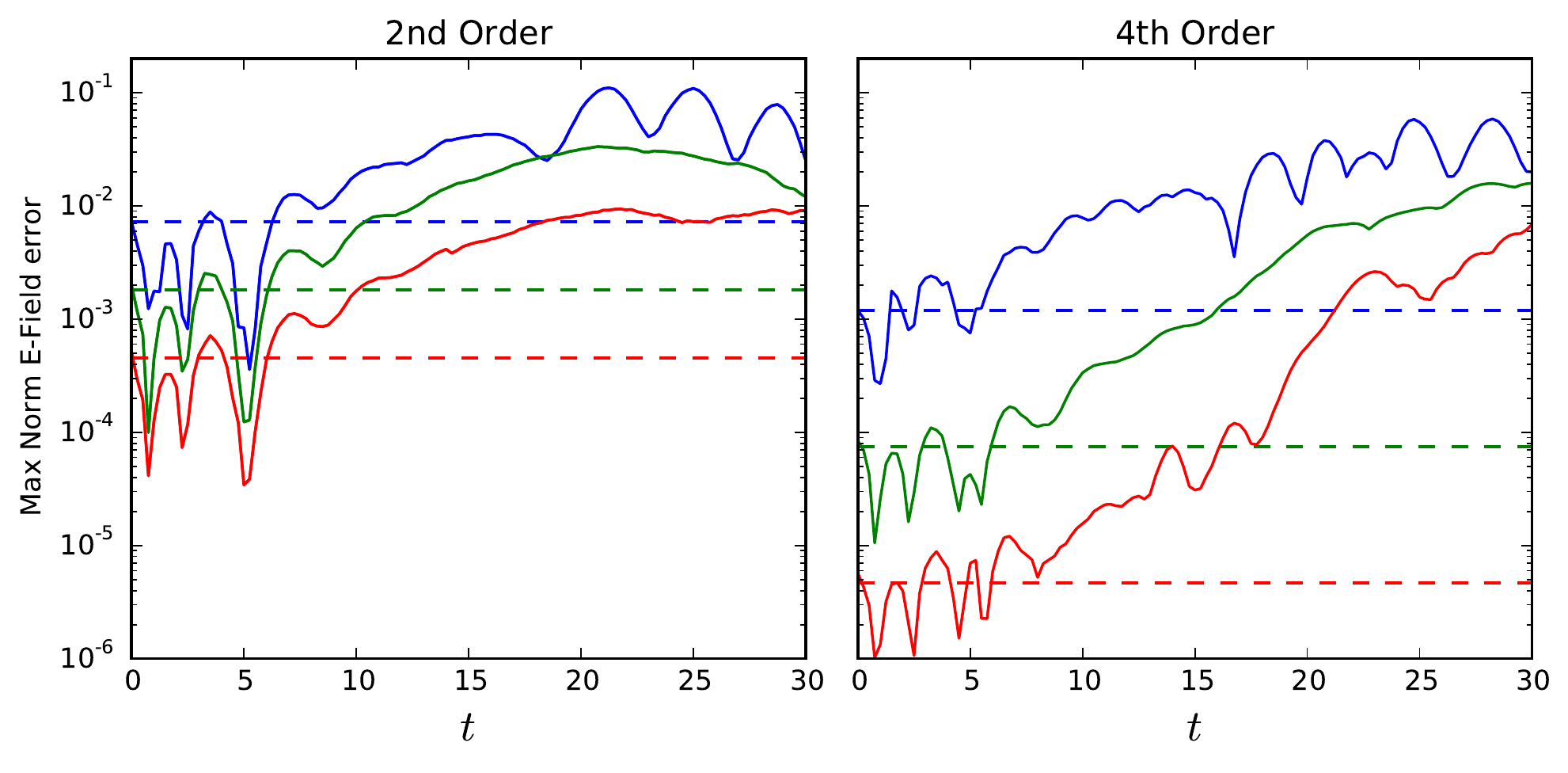}
    \caption{Max norm of the error in the electric field as a function of time for the 4th-order method on the two-dimensional, two-stream instability problem. Left - the 2nd-order method. Right - the 4th-order method. The line styles and colors have the same meaning as in Figure \ref{fig:landau_error}. \label{fig:twostream_error_2D}}
\end{figure}

Finally, we perform another convergence study, this time on the 2-dimensional setup. As before, we conduct 4 runs in total, starting at $N_{\text{cells}} = 8$, $N_x = 16$, and $N_v = 32$, and $dt = 1/4$, and increasing the resolution by a factor of 2 with each run. The resulting errors are shown as functions of time in Figure \ref{fig:twostream_error_2D}. As in the one-dimensional version of this problem, our method achieves 4th-order accuracy until late times, when under-resolved filaments reduce the accuracy of the solution to 2nd-order. Figure \ref{fig:twostream_error_2D} also compares these errors to those made by the 2nd-order method on the same problem setup. As in the one-dimensional case, the 4th-order method is significantly more accurate, particularly at early times. Overall, our method does not appear to suffer from unexpected numerical instabilities when operating in greater than one spatial dimension.
 
\section{Conclusions and Future Research}
\label{sec:conclusions}

We have presented a 4th-order accurate Particle-in-Cell algorithm for solving the Vlasov-Poisson equation in the context of electrostatic plasmas, including a remapping step that controls particle noise by periodically re-initializing the particle distribution on a Cartesian grid in phase space. We have demonstrated the accuracy of our method by comparing its performance to that of a 2nd-order method on a set of one-dimensional test problems. We have also gauged our method's performance on a set of two-dimensional test problems, finding that the 4th-order convergence of our method is maintained. 

In future work, we will explore extending the current method to 3 spatial and 3 velocity dimensions. While algorithmically straightforward, working in six total dimensions is computationally challenging, particularly during the particle remapping step, which presently involves working with high-dimensional grids. We will address this issue in two ways. First, we will explore phase-space aware parallelization strategies that more naturally map to the work distribution of warm plasma calculations. In contrast, our current parallelization strategy is based on domain decomposition in physical space only, such that the work per process grows more rapidly than the number of parallel domains as the problem size increases. This clearly limits the degree to which our implementation can scale to the high process counts that will be critical in six phase-space dimensions. Second, we will explore grid-free methods of performing the remap step. This will prevent us from needing to allocate temporary six-dimensional grids, which will greatly reduce the memory requirements and communication costs associated with the remap step.

\bibliographystyle{unsrt}
\bibliography{4thorderpic}

\begin{thebibliography}{10}

\bibitem{yt_paper}
M.~J. {Turk}, B.~D. {Smith}, J.~S. {Oishi}, S.~{Skory}, S.~W. {Skillman},
  T.~{Abel}, and M.~L. {Norman}.
\newblock {yt: A Multi-code Analysis Toolkit for Astrophysical Simulation
  Data}.
\newblock {\em The Astrophysical Journal Supplement Series}, 192:9, January
  2011.

\bibitem{scipy}
Eric Jones, Travis Oliphant, Pearu Peterson, et~al.
\newblock {SciPy}: Open source scientific tools for {Python}, 2001--.
\newblock [Online; accessed 2015-11-11].

\bibitem{Ipython}
Fernando P\'erez and Brian~E. Granger.
\newblock {IP}ython: a system for interactive scientific computing.
\newblock {\em Computing in Science and Engineering}, 9(3):21--29, May 2007.

\bibitem{NumPy}
S.~van~der Walt, S.C. Colbert, and G.~Varoquaux.
\newblock The numpy array: A structure for efficient numerical computation.
\newblock {\em Computing in Science and Engineering}, 13(2):22--30, March 2011.

\bibitem{Matplotlib}
J.D. Hunter.
\newblock Matplotlib: A 2d graphics environment.
\newblock {\em Computing in Science and Engineering}, 9(3):90--95, June 2007.

\bibitem{filbet_comparison_2003}
F.~{Filbet} and E.~{Sonnendr{\"u}cker}.
\newblock {Comparison of Eulerian Vlasov solvers}.
\newblock {\em Computer Physics Communications}, 150:247--266, February 2003.

\bibitem{banks_new_2010}
J.~W. {Banks} and J.~A.~F. {Hittinger}.
\newblock {A New Class of Nonlinear Finite-Volume Methods for Vlasov
  Simulation}.
\newblock {\em IEEE Transactions on Plasma Science}, 38:2198--2207, September
  2010.

\bibitem{vogman_dory_2014}
G.~V. {Vogman}, P.~{Colella}, and U.~{Shumlak}.
\newblock {Dory-Guest-Harris instability as a benchmark for continuum kinetic
  Vlasov-Poisson simulations of magnetized plasmas}.
\newblock {\em Journal of Computational Physics}, 277:101--120, November 2014.

\bibitem{hockney_computer_1981}
R.~W. Hockney and J.~W. Eastwood.
\newblock {\em Computer Simulation Using Particles}.
\newblock 1981.

\bibitem{heitmann_robustness_2005}
Katrin Heitmann, Paul~M. Ricker, Michael~S. Warren, and Salman Habib.
\newblock Robustness of cosmological simulations. i. large-scale structure.
\newblock {\em The Astrophysical Journal Supplement Series}, 160:28--58,
  September 2005.

\bibitem{warp}
D.~P. {Grote}, A.~{Friedman}, J.-L. {Vay}, and I.~{Haber}.
\newblock {The WARP Code: Modeling High Intensity Ion Beams}.
\newblock In M.~{Leitner}, editor, {\em Electron Cyclotron Resonance Ion
  Sources}, volume 749 of {\em American Institute of Physics Conference
  Series}, pages 55--58, March 2005.

\bibitem{edwards_2012}
E.~{Edwards} and R.~{Bridson}.
\newblock {A high-order accurate particle-in-cell method}.
\newblock {\em International Journal for Numerical Methods in Engineering},
  90:1073--1088, June 2012.

\bibitem{jacobs_2006}
G.~B. {Jacobs} and J.~S. {Hesthaven}.
\newblock {High-order nodal discontinuous Galerkin particle-in-cell method on
  unstructured grids}.
\newblock {\em Journal of Computational Physics}, 214:96--121, May 2006.

\bibitem{wang_particle--cell_2011}
B.~Wang, G.~Miller, and P.~Colella.
\newblock A particle-in-cell method with adaptive phase-space remapping for
  kinetic plasmas.
\newblock {\em {SIAM} Journal on Scientific Computing}, 33(6):3509--3537,
  January 2011.

\bibitem{cottet_vortex_2000}
Georges-Henri Cottet and Petros~D. Koumoutsakos.
\newblock {\em Vortex Methods: Theory and Practice}.
\newblock Cambridge University Press, March 2000.

\bibitem{chaniotis_remeshed_2002}
A.~K. Chaniotis, D.~Poulikakos, and P.~Koumoutsakos.
\newblock Remeshed smoothed particle hydrodynamics for the simulation of
  viscous and heat conducting flows.
\newblock {\em Journal of Computational Physics}, 182:67--90, October 2002.

\bibitem{wang_adaptive_2012}
B.~Wang, G.~Miller, and P.~Colella.
\newblock An adaptive, high-order phase-space remapping for the two dimensional
  vlasov--poisson equations.
\newblock {\em {SIAM} Journal on Scientific Computing}, 34(6):B909--B924,
  January 2012.

\bibitem{myers2015}
A.~{Myers}, P.~{Colella}, and B.~{Van Straalen}.
\newblock {The Convergence of Particle-in-Cell Schemes for Cosmological Dark
  Matter Simulations}.
\newblock {\em ArXiv e-prints}, March 2015.

\bibitem{williams_adaptive_2009}
Samuel Williams, Andrew Waterman, and David Patterson.
\newblock Roofline: An insightful visual performance model for multicore
  architectures.
\newblock {\em Commun. ACM}, 52(4):65--76, April 2009.

\bibitem{lo_real_2015}
B.~{Lo}, V.~{Minden}, and P.~{Colella}.
\newblock {A Real-Space Green's Function Method for Numerical Solution of
  Maxwell's Equations}.
\newblock {\em in preparation}, 2015.

\bibitem{chombo_design}
M.~Adams, P.~Colella, D.~T. Graves, J.N. Johnson, N.D. Keen, T.~J. Ligocki,
  D.~F. Martin, P.W. McCorquodale, D.~Modiano, P.O. Schwartz, T.D. Sternberg,
  and B.~Van Straalen.
\newblock \textit{Chombo Software Package for AMR Applications - Design
  Document}, {Lawrence Berkeley National Laboratory Technical Report
  LBNL-6616E}.

\bibitem{monaghan_particle_1985}
J.~J. Monaghan.
\newblock Particle methods for hydrodynamics.
\newblock {\em Computer Physics Reports}, 3:71--124, October 1985.

\bibitem{abramowitz_handbook_1972}
M.~{Abramowitz} and I.~A. {Stegun}.
\newblock {\em {Handbook of Mathematical Functions}}.
\newblock 1972.

\bibitem{miniati_block_2007}
F.~{Miniati} and P.~{Colella}.
\newblock {Block structured adaptive mesh and time refinement for hybrid,
  hyperbolic + N-body systems}.
\newblock {\em Journal of Computational Physics}, 227:400--430, November 2007.

\end{thebibliography}
    
\end{document}